\documentclass{ifacconf}

\usepackage{graphicx}      
\usepackage{natbib}        

\makeatletter
\let\old@ssect\@ssect 
\makeatother

\usepackage{hyperref}
\hypersetup{hidelinks}

\makeatletter
\def\@ssect#1#2#3#4#5#6{%
  \NR@gettitle{#6}
  \old@ssect{#1}{#2}{#3}{#4}{#5}{#6}
}
\makeatother

\usepackage{amsmath}
\usepackage{amssymb}
\usepackage{amsfonts, latexsym, amscd, enumerate, bbm, nicefrac, mathrsfs}
\usepackage{optidef}
\usepackage{xcolor}

\usepackage{listings}

\definecolor{codegreen}{rgb}{0,0.6,0}
\definecolor{codegray}{rgb}{0.5,0.5,0.5}
\definecolor{codered}{rgb}{0.65,0.0,0.0}
\definecolor{codepurple}{rgb}{0.58,0,0.82}
\definecolor{backcolour}{rgb}{0.95,0.95,0.92}
\definecolor{shadecolor}{gray}{0.85}
\definecolor{KULlightblue}{HTML}{116E8A}
\definecolor{KULdarkblue}{HTML}{00407A}

\lstdefinestyle{mystyle}{
    backgroundcolor=\color{backcolour},   
    commentstyle=\color{codegreen},
    keywordstyle=\color{magenta},
    numberstyle=\tiny\color{codegray},
    stringstyle=\color{codepurple},
    basicstyle=\ttfamily\footnotesize,
    breakatwhitespace=false,         
    breaklines=true,                 
    captionpos=b,                    
    keepspaces=true,                 
    numbers=none,                    
    numbersep=5pt,                  
    showspaces=false,                
    showstringspaces=false,
    showtabs=false,                  
    tabsize=2
}

\lstdefinestyle{Cstyle}{language=C,
    commentstyle=\color{KULdarkblue},
    keywordstyle=\color{KULdarkblue},
    numberstyle=\tiny\color{codegray},
    stringstyle=\color{codered},
    basicstyle=\ttfamily\tiny,
    breakatwhitespace=false,         
    breaklines=true,                 
    captionpos=b,                    
    keepspaces=false,           
    columns=fullflexible,
    numbers=none,                    
    numbersep=5pt,                  
    showspaces=false,                
    showstringspaces=false,
    showtabs=false,                  
    tabsize=1,
    xleftmargin=.0\textwidth, 
    xrightmargin=.0\textwidth
}

\lstdefinestyle{pythonstyle}{language=python,
    commentstyle=\color{KULdarkblue},
    keywordstyle=\color{KULdarkblue},
    numberstyle=\tiny\color{codegray},
    stringstyle=\color{codered},
    basicstyle=\ttfamily\tiny,
    breakatwhitespace=false,         
    breaklines=true,                 
    captionpos=b,       
    keepspaces=true,
    columns=fullflexible,
    numbers=none,                    
    numbersep=5pt,                  
    showspaces=false,                
    showstringspaces=false,
    showtabs=false,   
    tabsize=2,
    xleftmargin=.00\textwidth, 
    xrightmargin=.00\textwidth,
}

\lstdefinestyle{Cinline}{language=C,
    commentstyle=\color{KULdarkblue},
    keywordstyle=\color{KULdarkblue},
    numberstyle=\tiny\color{codegray},
    stringstyle=\color{codered},
    basicstyle=\ttfamily\small,
    breakatwhitespace=false,         
    breaklines=true,                 
    captionpos=b,                    
    keepspaces=false,           
    columns=fullflexible,
    numbers=none,                    
    numbersep=5pt,                  
    showspaces=false,                
    showstringspaces=false,
    showtabs=false,                  
    tabsize=1,
    xleftmargin=.0\textwidth, 
    xrightmargin=.0\textwidth
}

\lstdefinestyle{pythoninline}{language=python,
    commentstyle=\color{KULdarkblue},
    keywordstyle=\color{KULdarkblue},
    numberstyle=\tiny\color{codegray},
    stringstyle=\color{codered},
    basicstyle=\ttfamily\small,
    breakatwhitespace=false,         
    breaklines=true,                 
    captionpos=b,       
    keepspaces=true,
    columns=fullflexible,
    numbers=none,                    
    numbersep=5pt,                  
    showspaces=false,                
    showstringspaces=false,
    showtabs=false,   
    tabsize=2,
    xleftmargin=.00\textwidth, 
    xrightmargin=.00\textwidth,
}

\newcommand\YAMLcolonstyle{\color{black}\mdseries}
\newcommand\YAMLkeystyle{\color{green}\bfseries}
\newcommand\YAMLvaluestyle{\color{blue}\mdseries}

\makeatletter

\newcommand\language@yaml{yaml}

\expandafter\expandafter\expandafter\lstdefinelanguage
\expandafter{\language@yaml}
{
  keywords={true,false,null,y,n},
  keywordstyle=\color{darkgray}\bfseries,
  basicstyle=\ttfamily\tiny,
  sensitive=false,
  numbersep=0pt,    
  comment=[l]{\#},
  morecomment=[s]{/*}{*/},
  commentstyle=\color{purple}\ttfamily,
  stringstyle=\YAMLvaluestyle\ttfamily,
  moredelim=[l][\color{orange}]{\&},
  moredelim=**[il][\YAMLcolonstyle{:}\YAMLvaluestyle]{:},   
  morestring=[b]',
  morestring=[b]",
  literate =    {---}{{\ProcessThreeDashes}}3
                {>}{{\textcolor{red}\textgreater}}1     
                {|}{{\textcolor{red}\textbar}}1 
                {\ -\ }{{\mdseries\ -\ }}3,
    xleftmargin=0\textwidth, 
    xrightmargin=0\textwidth
}

\lst@AddToHook{EveryLine}{\ifx\lst@language\language@yaml\YAMLkeystyle\fi}
\makeatother
\newcommand\ProcessThreeDashes{\llap{\color{cyan}\mdseries-{-}-}}

\newcommand{\codename}[1]{{\fontfamily{pcr}\selectfont #1}}
\newcommand*{\impact}{{\codename{IMPACT}} }
\newcommand*{\casadi}{{\codename{CasADi}} }

\begin{document}
\begin{frontmatter}

\title{IMPACT: A Toolchain for Nonlinear Model Predictive Control Specification, Prototyping, and Deployment 
\thanksref{footnoteinfo}} 

\thanks[footnoteinfo]{This work was supported by the project Flanders Make SBO DIRAC: ``Deterministic and Inexpensive Realizations of Advanced Control". 
Flanders Make is the Flemish strategic research centre for the manufacturing industry.
}

\author[First]{Alvaro Florez\thanksref{equal_contrib}}
\author[First]{Alejandro Astudillo\thanksref{equal_contrib}}
\author[First]{Wilm Decré}
\author[First]{Jan Swevers}
\author[First]{Joris Gillis}

\address[First]{MECO~Research~Team,~Dept.~of~Mechanical~Engineering,~KU~Leuven.\\
Flanders Make@KU Leuven, 3001 Leuven, Belgium.\\
    (e-mail: \{alvaro.florez, alejandro.astudillovigoya, wilm.decre, jan.swevers, joris.gillis\}@kuleuven.be).}

\thanks[equal_contrib]{A. Florez and A. Astudillo contributed equally to this work.}

\begin{abstract}    
We present \codename{IMPACT}, a flexible toolchain for nonlinear model predictive control (NMPC) specification with automatic code generation capabilities. The toolchain reduces the engineering complexity of NMPC implementations by providing the user with an easy-to-use application programming interface, and with the flexibility of using multiple state-of-the-art tools and numerical optimization solvers for rapid prototyping of NMPC solutions. \codename{IMPACT} is written in \codename{Python}, users can call it from \codename{Python} and \codename{MATLAB}, and the generated NMPC solvers can be directly executed from \codename{C}, \codename{Python}, \codename{MATLAB} and \codename{Simulink}. An application example is presented involving problem specification and deployment on embedded hardware using \codename{Simulink}, showing the effectiveness and applicability of \codename{IMPACT} for NMPC-based solutions.
\end{abstract}


\end{frontmatter}

\section{Introduction} \label{sec:introduction}

The burgeoning computational power of off-the-shelf processing units has enabled the use of model predictive control (MPC) for applications that involve increasingly complex systems, need to comply with constraints and need to account for several, sometimes conflicting, performance objectives. However, development and implementation of nonlinear MPC (NMPC) comes with a high engineering cost, and the availability of flexible, easy-to-use tools for rapid prototyping and deployment of NMPC solvers is insufficient. This has limited the wide adoption of NMPC, e.g., for complex and fast mechatronic applications.


Central to an MPC control scheme is the repeated online numerical solution of a parametric optimal control problem (OCP). Important exceptions such as \codename{GRAMPC} \citep{grampc} notwithstanding, most efficient solvers rely on a direct transcription method to transform the OCP into a nonlinear program (NLP).
In recent years, a rich landscape of efficient quadratic programming (QP) solver implementations has emerged. They are valuable as a subroutine in 
an important 
class of NLP solvers, as well as directly in linear MPC, and real-time iteration schemes for NMPC.

No combination of OCP method, NLP solver and QP solver is applicable, let alone optimal, across a wide range of control applications.
This situation has given rise to a rich solver landscape\footnote{See 
\url{https://github.com/meco-group/dynamic_optimization_inventory}
}.
Some of these solvers are akin to a framework: offering some degree of vertical integration, i.e., going all the way from high-level easy-to-use OCP modeling down to low-level efficient deployable code.
Examples of efficient frameworks are: \codename{AutoGenU} \citep{autogenu}, providing an environment to formulate, solve, simulate and code generate MPC within \codename{Jupyter} notebooks, featuring a generalized minimal residual method implemented in \codename{C++} --- \codename{OpEn} \citep{OpEn2020}, focusing on embedded nonconvex optimization using the proximal averaged Newton-type method for optimal control (\codename{PANOC}) and the augmented Lagrangian method, implemented in \codename{Rust} and providing interfaces to \codename{Python} and \codename{MATLAB} --- \codename{Acados} \citep{ACADOS} and \codename{MATMPC} \citep{MATMPC}, offering SQP-type methods with \codename{Simulink} integration --- \codename{ParNMPC} \citep{ParNMPC}, offering an interior point solver with parallelization capabilities in \codename{MATLAB} --- \codename{MATLAB} MPC toolbox \citep{MPCmatlab}, offering a graphical interface and tight integration with \codename{Simulink}.

The main drawback of existing frameworks is that they are monolithic rather than modular: they restrict the user to a single or limited set of numerical algorithms, conflicting with the need of practitioners to prototype, e.g. figuring out the best solver and settings combination for a specific OCP. Furthermore, each provides a specific way to model OCPs and to communicate with the generated solution, which makes switching between frameworks, for testing and prototyping, a cumbersome process. Lastly, a lack of high-level modeling means that the user is burdened with writing in canonical forms and supplying derivatives.

To address the aforementioned issues, we present \codename{IMPACT}, 
an open-source toolchain that aims to facilitate the workflow of specification, prototyping and deployment of NMPC. It offers vertical integration, while providing modularity, with minimal runtime overhead.
The main features of \impact are described as follows:
\begin{itemize}
\item It provides a unified high-level application programming interface (API) to model OCPs.
\item It is cross-lingual (written in \codename{Python} and provides auto-generated\footnote{See
\url{https://gitlab.kuleuven.be/meco-software/python_matlab}
} \codename{MATLAB} bindings) and cross-platform.
\item It offers a range of direct transcriptions methods that 
are
paired to state-of-the-art QP/NLP solvers via \casadi plugins such as \codename{HPIPM} \citep{HPIPM},
\codename{PROX-QP} \citep{ProxQP},
\codename{IPOPT} \citep{ipopt},
\codename{FSLP} \citep{Kiessling2022CDC}.
\item It offers alternatives to built-in direct transcription by interfacing OCP solvers/frameworks via \codename{Rockit} plugins -- at the moment \codename{Acados}, \codename{GRAMPC} and \codename{FATROP} \citep{FatropIcra}.
\item It exports problem-specific MPC libraries with a unified \codename{C} API, useful for rapid prototyping and deployment on simulated and real environments.
\item It exports a problem-specific \codename{Simulink} block making use of the exported MPC library.
\end{itemize}
\impact is released under the LGPLv3 license and its source code is available at \url{https://gitlab.kuleuven.be/meco-software/impact}. It is built on top of the powerful tools for optimization and optimal control \casadi \citep{casadi} and \codename{Rockit} \citep{rockit}.

The remainder of this paper is organized as follows. Section \ref{sec:preliminaries} presents preliminary concepts on the formulation and transcription of the OCP underpinning MPC. In Section \ref{sec:toolchain_overview}, details on the overall structure and workflow of \impact are given. An application example using \impact is described in Section \ref{sec:case}. Finally, we close the paper with conclusions and future work remarks.

\section{Optimal Control Problem} \label{sec:preliminaries}
For a time horizon $t \in [t_0,\ t_{\mathrm{f}}]$, system state vector $\mathbf{x}(t) \in \mathbb{R}^{n_{\mathbf{x}}}$, input vector $\mathbf{u}(t) \in \mathbb{R}^{n_{\mathbf{u}}}$, algebraic  state vector $\mathbf{z}(t) \in \mathbb{R}^{n_{\mathbf{z}}}$, and parameter vector $\mathbf{p} \in \mathbb{R}^{n_{\mathbf{p}}}$, in \impact we consider general OCP formulations  of the canonical form
\begin{mini!}[2]
{{{\mathbf{x}, \mathbf{u}, \mathbf{z}}}}
{\int_{t_0}^{t_\mathrm{f}}V(\mathbf{x}(t), \mathbf{u}(t), \mathbf{p}, t)\mathrm{d}t + V_{f}(\mathbf{x}(t_\mathrm{f}),\mathbf{p}) \label{eq:ocp_objective}}{\label{eq:ocp}}{}
\addConstraint{B(\mathbf{x}(t_0),\mathbf{x}(t_f),\mathbf{p})}{ \leq 0 \label{eq:ocp_init}}{}
\addConstraint{\mathbf{\dot{x}}(t) }{ = \xi(\mathbf{x}(t),\mathbf{u}(t), \mathbf{z}(t),\mathbf{p}),\hspace{1.0ex} t \in [t_0, t_\mathrm{f}] \label{eq:ocp_dynamics}}{}
\addConstraint{\Gamma(\mathbf{x}(t),\mathbf{u}(t), \mathbf{z}(t),\mathbf{p}) }{ = 0,\hspace{3.4ex} t \in [t_0, t_\mathrm{f}] \label{eq:ocp_algebraic}}{}
\addConstraint{h(\mathbf{x}(t), \mathbf{u}(t), \mathbf{p})}{\leq 0,\hspace{8.2ex} t\in [t_0, t_\mathrm{f}], \label{eq:ocp_path}}{}
\end{mini!}
where $V(\cdot)$ and $V_\mathrm{f}(\cdot)$ in \eqref{eq:ocp_objective} are smooth nonlinear functionals that define Lagrange and Mayer terms, respectively, $\mathbf{p}$ is a parameter vector that typically contains a vector of state measurements or estimations $\mathbf{x}_{\mathrm{meas}} \in \mathbb{R}^{n_{\mathbf{x}}}$, \eqref{eq:ocp_init} defines a boundary constraint that typically takes the form of $\mathbf{x}(t_0)=\mathbf{x}_{\mathrm{meas}}$, \eqref{eq:ocp_dynamics} and \eqref{eq:ocp_algebraic} define a system of differential-algebraic equations (DAE) representing the system model, while \eqref{eq:ocp_path} represent general path constraints. Note that the inclusion of the algebraic state vector $\mathbf{z}$ and the algebraic equation \eqref{eq:ocp_algebraic} is not required and can be omitted for systems represented solely by ordinary differential equations (ODE). The OCP specified by the user in \impact (see Section \ref{sec:workflow}) is automatically transformed into the canonical form \eqref{eq:ocp}.

Direct transcription of \eqref{eq:ocp} leads to a finite-dimensional optimization problem
of the form
\begin{equation}
\label{eq:nlp}
\begin{split}
\min_{\mathbf{w}} \: f(\mathbf{w})\quad\mathrm{s.t.} \:\: g(\mathbf{w})=0,\:\: h(\mathbf{w})\leq 0,
\end{split}
\end{equation}
where 
$\mathbf{w} \in \mathbb{R}^{n_{\mathbf{w}}}$ is the vector of decision variables, 
and $N := T/\delta_t \in \mathbb{N}$ is the number of discretization points within the prediction horizon $T := (t_{\mathrm{f}} - t_0) \in \mathbb{R}_{>0}$ with a sampling time $\delta_t \in \mathbb{R}_{>0}$. The first- and second-order derivatives arising in the optimality conditions of  \eqref{eq:nlp} are typically large and sparse.

\section{Toolchain Description} \label{sec:toolchain_overview}

This section presents an overview of the structure
of the \impact toolchain. First, a description of the fundamental dependencies of the toolchain is given. Next, the elements of the toolchain, and a workflow from MPC specification to deployment are presented. Finally, 
the use of 
the unified API to communicate with the MPC solvers is detailed.

\subsection{Fundamental Dependencies}

\label{sec:dependencies}

%
The fundamental software dependencies required for installing and using \impact are: \codename{CasADi}, \codename{Rockit}, and the \codename{Python} modules \codename{pyyaml} and \codename{lxml}.

\impact relies on \codename{CasADi}
to handle 
expressions that model 
the objective and constraints as expression graphs, and to allow the saving/loading of functions and algorithms and their transfer across \codename{Python} and \codename{MATLAB}.

If instructed by the user, \impact relies on \codename{Rockit} to perform a chosen OCP transcription method. At the moment, multiple-shooting, single-shooting, direct collocation, and a B-spline method are supported.
Derivatives of the resulting NLP objective and constraints are obtained through \codename{CasADi}'s sparsity-exploiting algorithmic differentiation (AD), optionally made more efficient by \codename{C}-code generation.

\codename{Rockit} further provides interfaces to third-party OCP frameworks such as \codename{GRAMPC}, \codename{Acados} and \codename{FATROP}. Initially meant for offline usage, these interfaces are being refactored to use the \codename{CasADi} codegen \codename{C} API, such that they can be embedded in \codename{CasADi} expressions and are compatible with \codename{C}-code generation.
Lastly, \codename{Rockit} relies on \codename{CasADi} for importing functional mock-up units (FMU) for input/output representations of systems.

%
The libraries \codename{pyyaml} and \codename{lxml} are used within \impact to handle YAML and XML files required in the process of loading system models and generating \codename{Simulink} blocks when exporting the solution, respectively.

\subsection{Structure of the \codename{IMPACT} toolchain} \label{sec:structure}

\begin{figure}[htpb]
\centering
\includegraphics[width=0.33\textwidth]{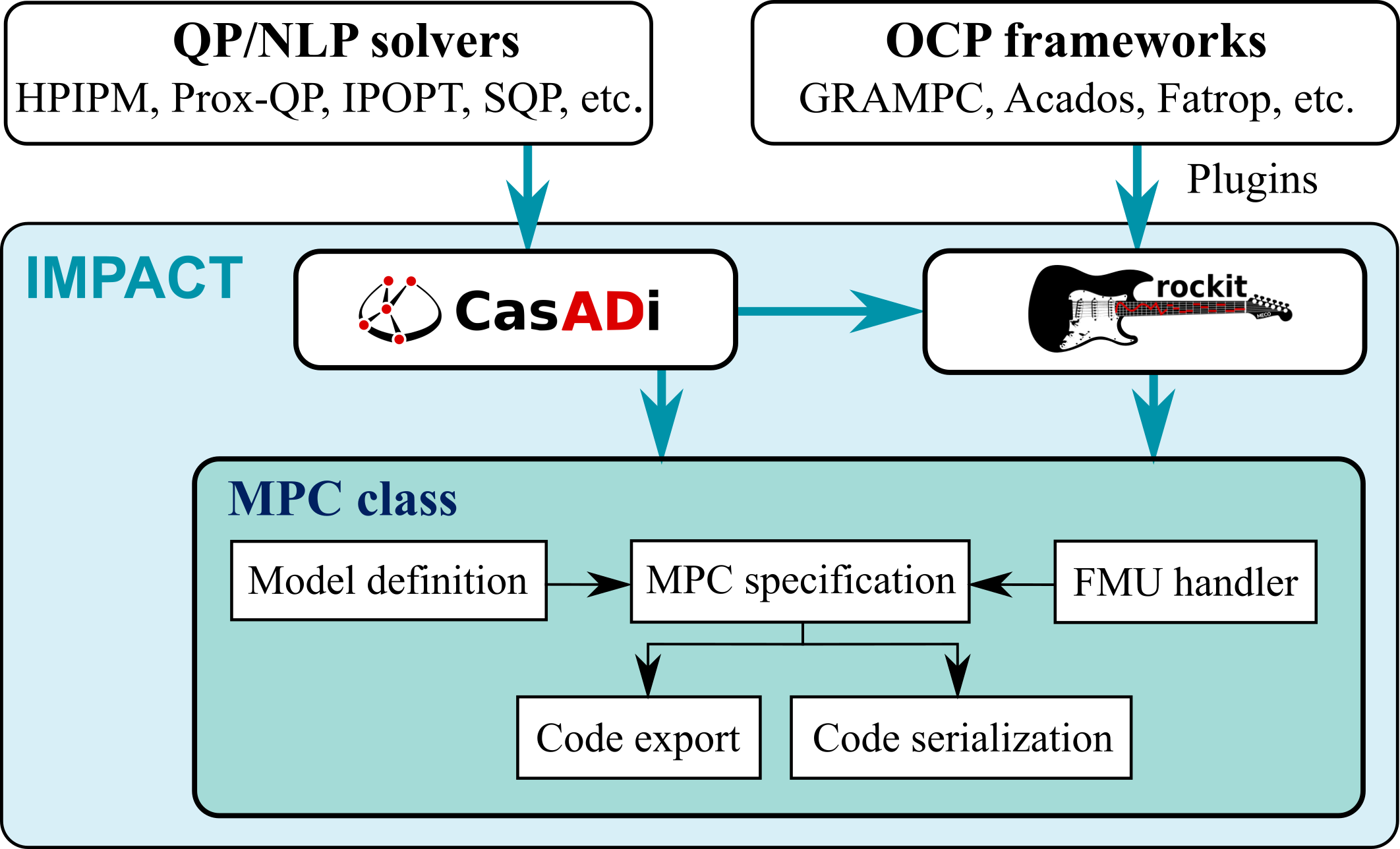}
\vspace{-1.5ex}
\caption{Overview architecture of the \impact toolchain showing the interaction between modules.}
\label{fig:structure}
\end{figure}

A graphical illustration of the overall structure of the \impact toolchain is shown in Fig. \ref{fig:structure}. 
%
\casadi and \codename{Rockit} play a vital role within \impact by providing essential capabilities to the toolchain regarding, e.g., expression handling,
and interfacing third-party solvers. In fact, the dependency of \impact on \casadi and \codename{Rockit}, and the inherent modularity that such dependencies represent, allow the direct and transparent availability of new solvers in \impact whenever an interface to such solvers is added to either \casadi or \codename{Rockit}.

Apart from 
\casadi and \codename{Rockit}, 
\impact is composed by the MPC class, which is instantiated by the user to use all the functionalities of \codename{IMPACT}. This class is divided into five sub-modules, which are described as follows.

\textit{1) Model definition:} \impact provides the user with two ways of defining (nonlinear) system models to be used within the OCP definition: (i) manually defining state and control variables in addition to differential and algebraic equations, as done in \codename{Rockit} or (ii) using YAML files -- i.e., which are human-readable files for data serialization -- by using the \lstinline[style=pythoninline]{add_model()} method of the MPC class. 
By using YAML files, we aim to standardize the definition of models based on ODEs or DAEs
for their use in any software tool that is compatible with YAML. 
This model representation allows the use of both inline definition of equations and external serialized \casadi functions, i.e., \lstinline[style=pythoninline]{*.casadi} files, as shown respectively in the two YAML snippets of a pendulum system that follow.
%

\vspace{-4ex}
\begin{minipage}[t]{.54\linewidth}
     \begin{lstlisting}[language=yaml]
equations:
  inline:
    ode:
      phi: dphi
      dphi: L*cos(phi)*sin(phi)...
differential_states: 
  - name: phi
  - name: dphi
controls: 
  - name: F
constants:
  inline:
    L: 2
            \end{lstlisting}
\end{minipage}\hfill
\begin{minipage}[t]{.44\linewidth}
      \begin{lstlisting}[language=yaml]
equations:
  external:
    type: casadi_serialized
    file_name: ode.casadi
differential_states: 
  - name: phi
  - name: dphi
controls: 
  - name: F
            \end{lstlisting}
\end{minipage}

\vspace{-2ex}

\textit{2) FMU handler:} This module allows the user to load FMUs -- an industry standard to define containers that represent black-box (input/output) models -- as \casadi functions to be used within the OCP definition. By calling the \lstinline[style=pythoninline]{add_simulink_fmu()} method of the MPC class, the FMU is wrapped into a \codename{C++} file, which is then compiled and loaded back as a \casadi external function. When the FMU provides derivatives (forward derivatives per FMI standard 2.0), these can be used. Otherwise, \casadi uses finite differences. In a further stage, FMU handling should be absorbed into the YAML file approach.

\textit{3) MPC specification:} This module inherits the functionality of the OCP class of \codename{Rockit} to enable the specification of the OCPs underpinning MPC. With the OCP class, \codename{Rockit} allows the definition of (multi-stage) OCPs by (i) managing symbolic variables for system states $\mathbf{x}$, system inputs $\mathbf{u}$, and algebraic states $\mathbf{z}$, (ii) associating dynamics to $\mathbf{x}$ and algebraic equations to $\mathbf{z}$ as in \eqref{eq:ocp_dynamics} and \eqref{eq:ocp_algebraic}, (iii) composing the objective function \eqref{eq:ocp_objective} by adding functions evaluated at different instants within the OCP horizon, and (iv) setting path or boundary constraints to expressions or variables as in \eqref{eq:ocp_init} and \eqref{eq:ocp_path}. When choosing a third-party OCP plugin, restrictions on the problem specification may apply.
%
The MPC specification is also used for result post-processing capabilities which allow the user, for instance, to retrieve the values of specific variables or expressions in the OCP solution, or to interpolate the solution on a refined grid of integration points.

\textit{4) Code serialization:} This module allows the serialization of the expression graph created by \casadi to define the MPC solver. This allows to save an MPC solver instantiation (including the functions that define the OCP) into a 
file that can be then loaded from \codename{Python}, \codename{MATLAB} or \codename{C++}. This module relies on the \lstinline[style=pythoninline]{save()} and \lstinline[style=pythoninline]{load()} methods of the MPC class.
In a similar fashion, instances of the \codename{IMPACT} Python module can also be saved and loaded.

\textit{5) Code export:} 
By calling the \lstinline[style=pythoninline]{export()} method of the MPC class, the user generates multiple artifacts that allow the prototyping and deployment of a solver within an MPC implementation.
The main artifacts depicted in Figure \ref{fig:workflow} are (i) an \impact library callable through the \impact \codename{C} API, (ii) a numerical backend library containing a concrete low-level problem formulation coupled to a solver that is either included as fully self-contained \codename{C} code, linked in with a dependency on the \codename{C++} \casadi runtime library and chosen solver plugin such as \codename{IPOPT}, or linked in as third-party dependencies such as \codename{OSQP} or \codename{GRAMPC}, and (iii) an MPC-ready \impact \codename{Simulink} block. 
%
%

The \impact library is accompanied by \codename{Python} bindings. More details can be found in Section \ref{sec:C_API}.

Self-contained \codename{C} code
is dependency-free and requires only a \codename{C} compiler to be deployed in any compatible target. It requires every element of the MPC solver instantiation to be compatible with code-generation -- e.g., if the user selects \codename{IPOPT}, which cannot be code-generated, no self-contained \codename{C} code of the MPC solver can be generated.

The 
\impact \codename{Simulink} block is based on a code-generated custom \codename{C} S-function. \impact generates a \lstinline[style=pythoninline]{.slx} container based on the \codename{C} API  and linked to the generated \codename{C} files, and can be directly used within a \codename{Simulink} model. \codename{Simulink} has become a de-facto tool for simulation and prototyping of control systems. This artifact allows \impact users to easily prototype their designed MPC with no extra steps. Moreover, if the linked \codename{C} file is self-contained, the generated \codename{Simulink} block would be compatible with the \codename{Simulink} coder. This way, the MPC-ready \impact \codename{Simulink} block could be code-generated as part of a larger \codename{Simulink} model and deployed into real-time targets without additional engineering effort.
Compatibility with \codename{Simulink} coder is also expected when third-party dependencies are open-source and \codename{C} based, but this worflow has not yet been validated.

The three artifacts include \casadi functions for the MPC solver (with parametric and hot-start inputs), forward simulation of the system dynamics, DAEs, performance objective, among others. In addition, they output statistics of the solution such as solution time, status, optimality of the solution and number of iterations. These provide the user with a great deal of flexibility to implement and prototype the MPC solution within a simulated or real control loop. Note that once the artifacts are generated, the only values that can be modified within the problem are the parameters $\mathbf{p}$.




%
%

\subsection{Toolchain Workflow}
\label{sec:workflow}
\begin{figure}[htpb]
\centering
\includegraphics[width=0.43\textwidth]{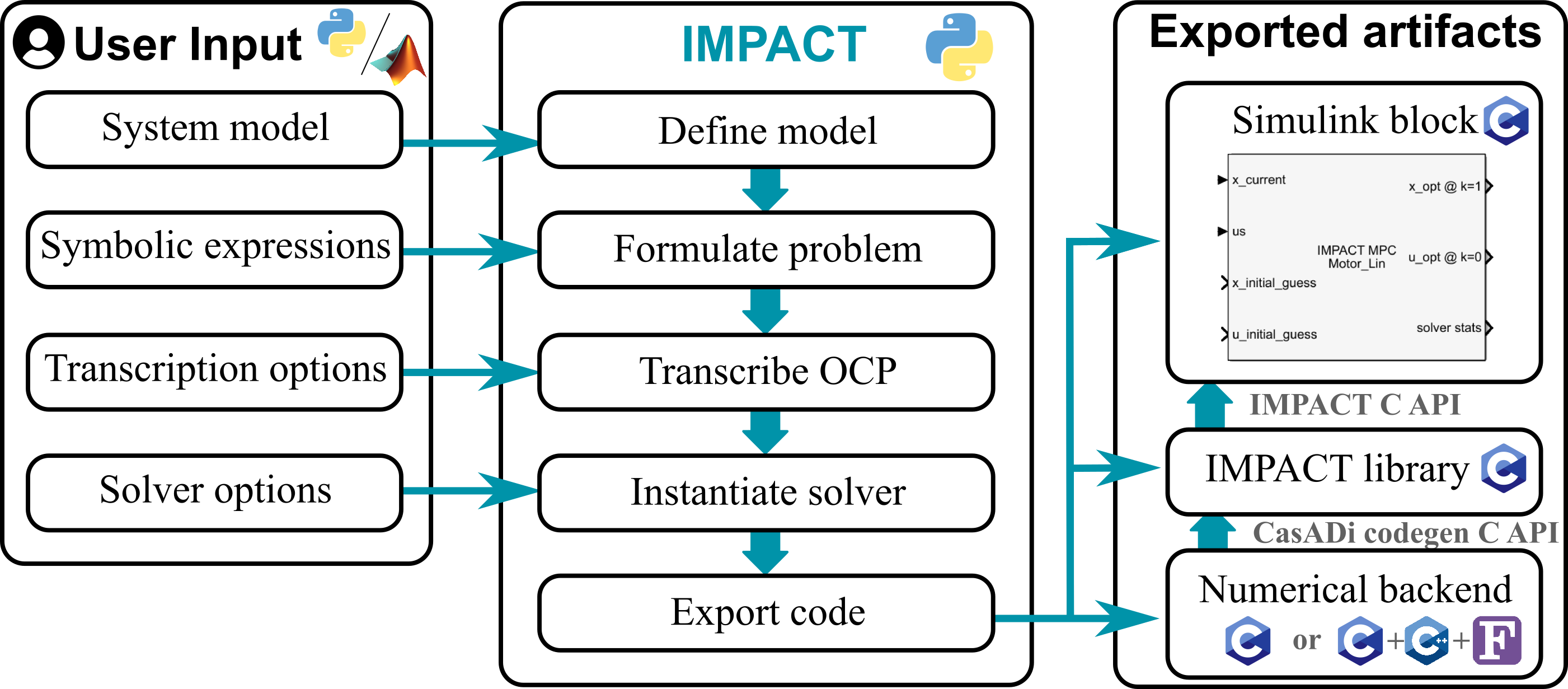}
\vspace{-1.5ex}
\caption{Overview of the workflow of the \impact toolchain.}
\label{fig:workflow}
\end{figure}

This section has given details on the architecture of \impact and its dependencies. It is now necessary to explain the workflow proposed within \impact for the user to specify an MPC, prototype its solution, and deploy it. A general overview of the workflow is presented in Fig. \ref{fig:workflow}.

The workflow in the \impact toolbox is described as follows by means of a simple example involving the control of a pendulum. A description of each step is provided followed by the code snippet associated with it. 

First of all, the user must import the \impact module and instantiate an object of the MPC class. The argument \lstinline[style=pythoninline]{T=2.0} sets the prediction horizon $T$ to a fixed value. However, $T$ could be considered as a decision variable within a free time problem by passing \lstinline[style=pythoninline]{T=FreeTime(2.0)} as argument, where \lstinline[style=pythoninline]{2.0} is an initial guess for this variable.
{
\vspace{-0.5ex}
\begin{lstlisting}[style=pythonstyle]
from impact import *
mpc = MPC(T=2.0)
\end{lstlisting}
\vspace{-1ex}
}
Afterwards, the system model should be declared. As mentioned in Section \ref{sec:structure}, this can be done by manually writing the ODE or DAE -- i.e., using instances of states and inputs within the MPC class --, or by loading a YAML file, which is the option shown in the code snippet below.
{
\vspace{-2ex}
\begin{lstlisting}[style=pythonstyle]
pendulum = mpc.add_model('pendulum','pendulum.yaml')
\end{lstlisting}
\vspace{-1.5ex}
}

Once the model has been defined, the user can instantiate the elements of the vector of parameters $\mathbf{p}$, such as the current state (\lstinline[style=pythoninline]{x_0}), a desired final state (\lstinline[style=pythoninline]{x_f}), and the weights (\lstinline[style=pythoninline]{Wt}) to be used in the objective.
{
\vspace{-0.5ex}
\begin{lstlisting}[style=pythonstyle]
x_0 = mpc.parameter('x_0',pendulum.nx)
x_f = mpc.parameter('x_f',pendulum.nx)
Wt  = mpc.parameter('Wt',2)
\end{lstlisting}
\vspace{-1.5ex}
}

In a high-level way, the objective and constraints are then defined 
by using the methods \lstinline[style=pythoninline]{mpc.add_objective()} and \lstinline[style=pythoninline]{mpc.subject_to()} directly using state and control names. Note that, by using \lstinline[style=pythoninline]{mpc.at_t0()} or \lstinline[style=pythoninline]{mpc.at_tf()} the user can specify boundary constraints.
\begin{lstlisting}[style=pythonstyle]
# Objective 
mpc.add_objective(mpc.integral(Wt[0]*pendulum.F**2 + Wt[1]*pendulum.dphi**2))
# Boundary constraints
mpc.subject_to(mpc.at_t0(pendulum.x) == x_0)
mpc.subject_to(mpc.at_tf(pendulum.x) == x_f)
# Path constraints
mpc.subject_to(-2 <= (pendulum.F <= 2 ))
\end{lstlisting}
\vspace{-1.5ex}
The nonlinear optimization solver to be used is defined by means of the \lstinline[style=pythoninline]{mpc.solver()} method, for which the user provides several 
options depending on the selected solver. In this case, the SQP method of \casadi has been selected.
\begin{lstlisting}[style=pythonstyle]
mpc.solver('sqpmethod',options={...})
\end{lstlisting}
\vspace{-1.5ex}
The transcription method is defined either by using a direct method provided by \codename{Rockit} or an external method provided by a plugin to an optimization tool such as \codename{GRAMPC}. Here, \lstinline[style=pythoninline]{N} corresponds to the number of discretization points $N$, while \lstinline[style=pythoninline]{intg='rk'} sets a $4$th-order Runge Kutta integrator to discretize the system dynamics.
\begin{lstlisting}[style=pythonstyle]
# Direct method, option 1
method = MultipleShooting(N=40,intg='rk')
# External method, option 2
method = external_method('grampc',N=40,grampc_options=...)
# Actually set the selected method within the MPC environment
mpc.method(method)
\end{lstlisting}
\vspace{-1.5ex}
The user can now set the values of the parameters $\mathbf{p}$ and execute the solver directly in \codename{Python}.
\begin{lstlisting}[style=pythonstyle]
mpc.set_value(x_0, [0.5,0,0,0]) # Set parameters
mpc.set_value(x_f, [0,0,0,0])
mpc.set_value(Wf, [1,1]))
solution = mpc.solve() # Solve the OCP
\end{lstlisting}
Finally, the user can generate the three artifacts explained in Section \ref{sec:structure} by executing the following code.
\begin{lstlisting}[style=pythonstyle]
mpc.export("pendulum")
\end{lstlisting}
\vspace{-1.5ex}
For rapid prototyping of the exported code the user can: (i) run the generated \texttt{hello\_world} examples in \codename{C} or \codename{Python} which use the \codename{C} API, or (ii) use the generated \codename{Simulink} block from the generated \lstinline[style=pythoninline]{.slx} file in a \codename{Simulink} model. Since we are using the SQP method of \codename{CasADi}, which is compatible with code-generation, in this example we could deploy the code-generated code or the \codename{Simulink} block in a real-time target architecture.


\subsection{\codename{IMPACT} \codename{C} API}
\label{sec:C_API}

The \codename{IMPACT} \codename{C} API defines a unified API to use the MPC solver exported by \codename{IMPACT}. It represents the core artifact of the code export functionality and is a layer of abstraction that makes the implementation agnostic to the content of the generated MPC solver, e.g., self-contained \codename{C} code with \codename{C} solver, \codename{C++} code that still depends on \codename{CasADi} and/or third-party solvers. 
It provides rich functions to set the inputs of the MPC solver (parameters) and get the values of variables from the optimized solution.
For the sake of completeness, we present below a minimum example of the use of the \impact \codename{C} API 
for the pendulum example presented in Section \ref{sec:workflow}.
\begin{lstlisting}[style=Cstyle]
impact_struct* m = impact_initialize();
impact_set(m, "x_0", IMPACT_ALL, 0, x_meas, IMPACT_FULL); // Set a parameter
impact_solve(m); // Solve a single OCP
// Get optimal input for first time instance
impact_get(m, "u_opt", IMPACT_ALL, 0, U, IMPACT_FULL);
\end{lstlisting}
\vspace{-1.5ex}
Here, \lstinline[style=Cinline]{impact_initialize()} instantiates the exported \impact library and returns a pointer \lstinline[style=Cinline]{m} to it.
By using the \lstinline[style=Cinline]{impact_set()} method, the parameter \lstinline[style=Cinline]{x_0} is set with a user-defined value \lstinline[style=Cinline]{x_meas}. 
The third argument \lstinline[style=Cinline]{IMPACT_ALL} tells \impact that all the elements of the specified parameter are set. The user can use \lstinline[style=Cinline]{"pendulum.phi"} instead of \lstinline[style=Cinline]{IMPACT_ALL} to set the parameter corresponding to the angle \lstinline[style=Cinline]{phi}, for instance.
The fourth argument \lstinline[style=Cinline]{0} defines the time instance within the horizon for which the parameter is set. 
Another option is \lstinline[style=Cinline]{IMPACT_EVERYWHERE}, which sets the parameter for the whole horizon. The sixth argument combines flags representing repetition (\lstinline[style=Cinline]{IMPACT_FULL}, \lstinline[style=Cinline]{IMPACT_HREP}) and data ordering (\lstinline[style=Cinline]{IMPACT_ROW_MINOR},  \lstinline[style=Cinline]{IMPACT_COLUMN_MAJOR}).
The \lstinline[style=Cinline]{impact_solve()} function executes the MPC solver. Then, an element of the solution -- i.e., the optimal control input \lstinline[style=Cinline]{u_opt} -- is retrieved and assigned to a variable \lstinline[style=Cinline]{U} by using \lstinline[style=Cinline]{impact_get()}, with the arguments mirroring \lstinline[style=Cinline]{impact_set()}. 
This example shows the solution of one OCP. 
To implement (N)MPC, the \lstinline[style=Cinline]{impact_set()}, \lstinline[style=Cinline]{impact_solve()}, and \lstinline[style=Cinline]{impact_get()} methods are executed in a loop. More methods are defined in the API, e.g., to retrieve statistics and perform debugging.

\section{Application Example Using IMPACT} \label{sec:case}

This section demonstrates the use and deployment of MPC using \codename{IMPACT}. We control the point-to-point motion (angular position) of a DC motor using a Speedgoat SN7233 real-time target machine. The control scheme estimates a constant disturbance that represents unmodeled dynamics -- e.g., input disturbances and static friction --, counteracts its effect, and achieves offset-free positioning.

\subsection{System Description}
The hardware has two components: the DC motor and the real-time target machine that executes the MPC algorithm and interfaces with the motor driver and a rotary encoder. 

The shaft of the DC motor is connected to a load and the rotary encoder, which provides a measurement $\mathbf{y}$ of the angular position of the rotor $\theta \in \mathbb{R}$. The motor has a driver that receives an analog voltage $v_{\mathrm{motor}} \in [-10,10]$
from the digital-to-analog converter (DAC) of the Speedgoat. 
For a state vector $\mathbf{x} := [\theta, \dot{\theta}]^T \in \mathbb{R}^2$ and a control input $\mathbf{u} := v_{\mathrm{motor}}$, a second order linear model for the motor 
\begin{equation} \label{eq:dynamics_motor}
    \dot{\mathbf{x}} = A\mathbf{x} + B\mathbf{u}
\end{equation}
has been experimentally identified using the \codename{LCtoolbox} \citep{Maarten2018LCToolbox}. This model contains one integrator and one real pole which models the linear motor friction, but does not consider the static friction component. 

\begin{figure}[htpb]
    \centering
    \includegraphics[width=0.36\textwidth]{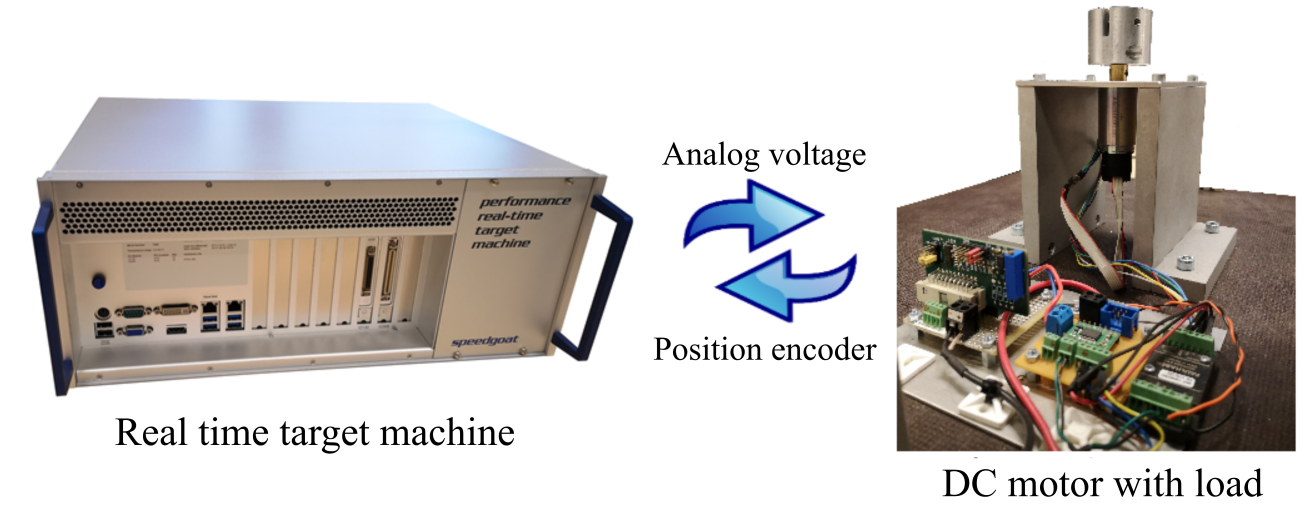}
    \vspace{-3ex}
    \caption{
    Overview of the system.
    }
    \label{fig:connectionScheme}
\end{figure}
The Speedgoat SN7233 is a real-time target machine, programmed via \codename{Simulink}, that allows the deployment of controllers on, e.g., multi-core CPUs, FPGAs and PLCs, and facilitates rapid controller prototyping and hardware-in-the-loop simulations through a smooth integration with \codename{Simulink Real-Time}. 
The target runs the MPC algorithm -- based on an \impact \codename{Simulink} block --, applies an analog voltage $v_{\mathrm{motor}}$ to the motor, and reads the digital signals from the rotary encoder, as depicted in Fig. \ref{fig:connectionScheme}. 

\subsection{Problem Definition} \label{sec:problem_definition}
\begin{figure}[htpb]
    \centering
    \includegraphics[width=0.37\textwidth]{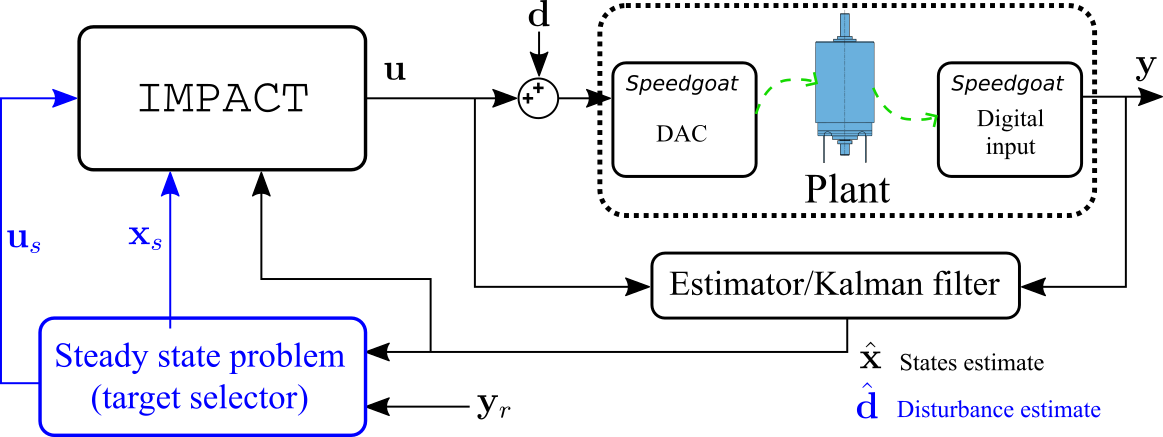}
    \vspace{-1.5ex}
    \caption{Diagram of an offset-free MPC implementation, with differences with respect to traditional MPC in blue, and physical connections to the motor in green.
    }
    \label{fig:offsetfree-diagram}
\end{figure}

To achieve an accurate positioning of the rotor, i.e., counteracting the effect of the model-plant mismatch and a constant input disturbance $\mathbf{d} \in \mathbb{R}$, we define an offset-free MPC \citep[Section 1.5]{rawlings2017model}.
In addition to the solution of OCP \eqref{eq:ocp},
this methodology involves (i) augmenting the system dynamics within the estimation algorithm with the dynamics of $\mathbf{d}$, and (ii) solving an optimization problem -- the steady-state problem or target selector -- that sets references for $\{\mathbf{x},\mathbf{u}\}$ in the OCP \eqref{eq:ocp} based on the estimated states and disturbance, to compensate the disturbance effects while tracking a reference. 
The offset-free MPC structure is synthesized in Fig. \ref{fig:offsetfree-diagram}. 

Before presenting details on the definition of OCP \eqref{eq:ocp} for this application, let us describe the estimator and steady-state problem. The estimator -- i.e., a Kalman filter -- outputs an estimate
$\hat{\bar{\mathbf{x}}}_k$ 
of an augmented state vector $\bar{\mathbf{x}}_k := \begin{bmatrix}\mathbf{x}_k^\top & \mathbf{d}_k \end{bmatrix}^\top$ whose dynamics are based on a discretized representation of \eqref{eq:dynamics_motor} and the assumption of a constant disturbance. 
The steady-state problem corresponds to an optimization problem that considers the reference $\mathbf{y}_r$ and the estimated disturbance $\hat{\mathbf{d}}$ to define set-points $\{\mathbf{x}_s,\mathbf{u}_s\}$ for $\{\mathbf{x},\mathbf{u}\}$ in OCP \eqref{eq:ocp}.
For this example, the solution of such problem can be found analytically as the solution of a system of linear equations. %
%
The reader is referred to 
\cite{rawlings2017model} for more details on this formulation.

Instead of minimizing the output error $\mathbf{y} - \mathbf{y}_r$, OCP \eqref{eq:ocp} aims to minimize the deviation between the 
set-points
$\{\mathbf{x}_s, \mathbf{u}_s\}$ and the system variables $\{\mathbf{x}, \mathbf{u}\}$. Therefore, the objective \eqref{eq:ocp_objective} is defined by
$V(\cdot) := \lVert\mathbf{x}(t)-\mathbf{x}_s\rVert_{Q}^2+\lVert\mathbf{u}(t)-\mathbf{u}_s\rVert_{R}^2$ and $V_{\mathrm{f}}(\cdot) := \lVert\mathbf{x}(t_f)-\mathbf{x}_s\rVert_{Q_\mathrm{f}}^2$, where $Q$, $R$, $Q_\mathrm{f} \succeq 0$ are weight matrices.
%
In the boundary constraint \eqref{eq:ocp_init}, $\mathbf{x}_{\mathrm{meas}} = \hat{\mathbf{x}}$,
i.e., excluding the estimated disturbance.
The dynamics \eqref{eq:ocp_dynamics} are described by the ODE \eqref{eq:dynamics_motor}. Since the system does not feature algebraic equations, \eqref{eq:ocp_algebraic} is omitted.
Boundary constraints \eqref{eq:ocp_path} are given by the lower and upper bounds on $v_{\mathrm{motor}}$, such that $-10 \leq \mathbf{u}(t) \leq 10\ \forall t \in [t_0, t_\mathrm{f}]$.



\subsection{Implementation Details}

To implement the \impact \codename{Simulink} block that solves OCP \eqref{eq:ocp} we follow the workflow outlined in Section \ref{sec:workflow}. We use $N = 50$,
$\delta_t = 3.33\ [\mathrm{ms}]$ and 
$T = 166.66\ [\mathrm{ms}]$.
Model \eqref{eq:dynamics_motor} is defined in a YAML file and loaded to the MPC object. Five parameters are defined, namely, $\mathbf{x}_{\mathrm{meas}}$, $\mathbf{x}_s$, $\mathbf{u}_s$, $Q$, $R$ and $Q_{\mathrm{f}}$. The objective and constraints of OCP \eqref{eq:ocp} are defined as in Section \ref{sec:problem_definition}. 
We use 
the active set-based QP solver \codename{QRQP} \citep{qrqp} in \codename{CasADi}
and multiple shooting as transcription method with a 4th-order Runge-Kutta integrator for system discretization. The MPC solver is exported to generate the 
\codename{Simulink} block. This block is loaded into a \codename{Simulink} model where the estimator and the steady-state problem are implemented, and where the Speedgoat input/output blocks are loaded. The \impact block is connected to the other blocks as in Fig. \ref{fig:offsetfree-diagram} to close the loop.
%
%
The reference is 
$\mathbf{y}_r(t) = -3\ [\mathrm{rad}]$, $\forall t \in [1, 3]$, and zero otherwise. 
While $\mathbf{d} = 0.5\ [\mathrm{V}]$ $\forall t \geq 8$.




\subsection{Results}
\begin{figure}[htpb]
    \centering
    \includegraphics[width=0.40\textwidth]{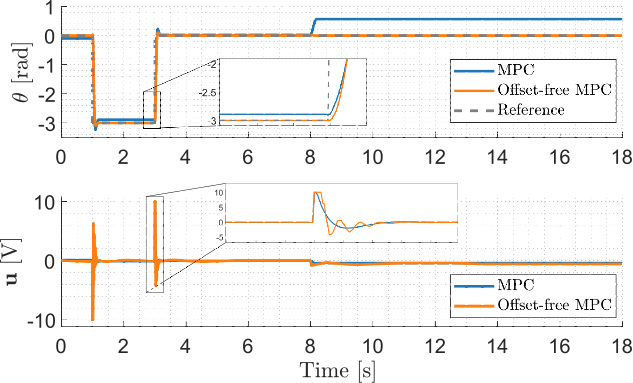}
    \vspace{-1.5ex}
    \caption{
    Evolution of the angular position $\theta$ and the voltage input $\mathbf{u}$ during the application execution.
    }
    \label{fig:1_posvolt2}
\end{figure}

Once the \codename{Simulink} model has been set, we call the real time execution operation of \codename{Simulink Real-Time}, which code-generates it by using the \codename{Simulink} coder and then executes it in the Speedgoat. The results of the offset-free MPC implementation are compared against a traditional MPC implementation, i.e., without disturbance estimation and execution of the steady-state problem.
%
%
Fig. \ref{fig:1_posvolt2} shows the evolution of the angular position $\theta$ and the voltage input $\mathbf{u}$ with both implementations. Here, the offset-free MPC can follow the reference with zero offset, while the traditional MPC deviates from the reference before (due to the effects of model-plant mismatch) and after the application of the disturbance at 
$t = 8$.


\impact eased MPC specification, testing and deployment, allowing to quickly prototype the problem and test several solvers by using the generated artifacts. In addition, the \impact \codename{Simulink} block could be easily integrated within larger models for simulation or deployment on real hardware.
%
%
Although this example features a linear system, \impact allows the definition of systems represented by both linear or nonlinear equations.
Other tested applications with \impact include the deployment of NMPC in a Speedgoat target and a Beckhoff TwinCAT Embedded PC for a point-to-point motion application with a parallel SCARA robot. In the future, \impact will be tested in a DSpace Controller board, and extended to be used in FPGAs and GPUs.

\section{Conclusions} \label{sec:conclusions}
In this paper we presented \codename{IMPACT}, a toolchain for nonlinear model predictive control that reduces the engineering complexity of the different steps from problem specification to MPC solution deployment. The different modules that compose the toolchain, their dependencies and their usage workflow were described. In the wide space of academic software, \impact stands out by offering a unified cross-lingual high-level way to model OCPs, and a unified \codename{C} API to communicate with exported MPC controllers. Future directions include the implementation of a method to allow the \codename{C} API to be exported as an FMU. We will also validate 
the toolchain in more complex scenarios involving complex robotic systems and self-driving cars.


{\small
\bibliography{impact} 
} 
                                                   







\end{document}